\documentclass[11pt]{article}
\usepackage{amscd}
\usepackage{amsmath}
\usepackage{amsfonts}

\newcommand{\dbar}{\overline{\partial}}
\newcommand{\dd}[1]{\partial_{#1}}
\newcommand{\dbr}[1]{\partial_{\overline{#1}}}

\newcommand{\ddt}[1]{\frac{\partial #1}{\partial t}}

\newcommand{\lt}{\tilde{\triangle}}
\newcommand{\chii}[2]{\chi^{#1 \overline{#2}}}

\newcommand{\chz}[2]{\chi_{0 \, #1 \overline{#2}}}
\newcommand{\g}[2]{g_{#1 \overline{#2}}}

\newcommand{\h}[2]{h^{#1 \overline{#2}}}
\renewcommand{\thesection}{1}

\newenvironment{proof}[1][Proof]{\begin{trivlist}
\item[\hskip \labelsep {\bfseries #1}]}{\end{trivlist}}

\begin{document}
\newcounter{theor}
\setcounter{theor}{1}
\newtheorem{claim}{Claim}
\newtheorem{theorem}{Theorem}[section]
\newtheorem{proposition}{Proposition}[section]
\newtheorem{lemma}{Lemma}[section]
\newtheorem{defn}{Definition}[theor]
\newtheorem{corollary}{Corollary}[section]

\centerline{\bf ON THE J-FLOW IN HIGHER DIMENSIONS }
\centerline{\bf AND THE LOWER BOUNDEDNESS OF}
\centerline{\bf THE MABUCHI ENERGY} 
\bigskip
\centerline{Ben Weinkove \footnote{This work was carried out while the author was
supported by a graduate fellowship at Columbia University.}}
\centerline{Department of Mathematics, Harvard University}
\centerline{1 Oxford Street, Cambridge, MA 02138}
\bigskip

\setlength\arraycolsep{2pt}
\addtocounter{section}{1}

\bigskip
\noindent
{\bf Abstract. } The J-flow is a parabolic flow on K\"ahler manifolds.  It was defined by Donaldson in the 
setting of moment maps and by Chen as the gradient flow of the J-functional appearing in his formula for the
Mabuchi energy.  It is shown here that under a certain condition on the initial data, the J-flow
converges to a critical metric. This is a generalization to higher dimensions of the author's previous work on
K\"ahler surfaces. A corollary of this is the lower boundedness of the Mabuchi energy on K\"ahler classes
satisfying a certain inequality when the first Chern class of the manifold is negative.

\bigskip
\noindent
{\bf 1. Introduction}
\bigskip

The J-flow is a parabolic flow of potentials on K\"ahler manifolds with two K\"ahler
classes.  It was defined by Donaldson
\cite{D1} in the setting of moment maps and by Chen \cite{C1} as the gradient
flow for the J-functional appearing in his
\cite{C1} formula for the Mabuchi energy \cite{Ma}.  Chen \cite{C2} proved
long-time existence of the flow for any smooth initial data, and proved
a convergence result in the case of non-negative bisectional curvature.  

The J-flow is defined as follows.  Let $(M, \omega)$ be a compact K\"ahler
manifold of complex dimension $n$ and let $\chi_0$ be another K\"ahler form on
$M$. Let $\mathcal{H}$ be the space of K\"ahler potentials
$$\mathcal{H} = \{ \phi \in C^{\infty}(M) \ | \ \chi_{\phi} = \chi_0 +
\frac{\sqrt{-1}}{2} \partial \dbar \phi >0 \}.$$
The J-flow is the flow in $\mathcal{H}$ given by
\begin{equation} \label{eqnJflow}
\left\{
\begin{array}{rcl}
{\displaystyle \ddt{\phi_t}} & = & \displaystyle{ c - \frac{\omega \wedge
\chi_{\phi_t}^{n-1}}{\chi_{\phi_t}^n}}
\\
\phi_0 & = & 0, \end{array} \right.
\end{equation}
where $c$ is the constant given by
$$c = \frac{ \int_M \omega \wedge \chi_0^{n-1}}{\int_M \chi_0^n}.$$
A critical point of the J-flow gives a K\"ahler metric $\chi$ satisfying
\begin{equation} \label{maineqn}
\omega \wedge \chi^{n-1} = c \chi^n.
\end{equation}
Donaldson \cite{D1} showed that 
a necessary condition for a solution to (\ref{maineqn}) in 
$[\chi_0]$ is that $[nc \chi_0 - \omega]$ be a positive class. 
He remarked that a natural conjecture would be that this condition be sufficient.  Chen \cite{C1}
proved this result if $n=2$, without using the J-flow.  
In \cite{W1}, the author gave an alternative proof by showing that for $n=2$, the
J-flow converges in 
$C^{\infty}$ to a critical metric
under the condition $nc \chi_0 - \omega>0$.

We generalize this as follows. 
 
\bigskip
\noindent
{\bf Main Theorem.} {\it If the K\"ahler metrics $\omega$ and $\chi_0$ satisfy
$$nc \chi_0 - (n-1)\omega >0,$$ 
then the J-flow (\ref{eqnJflow}) converges in $C^{\infty}$ to a smooth critical
metric.}

\bigskip 
This shows that (\ref{maineqn}) has a solution in $[\chi_0]$ under the 
condition $$nc [\chi_0] - (n-1) [\omega]>0.$$ 

Recall that the Mabuchi energy is a functional on $\mathcal{H}$ defined
by
$$M_{\chi_0}(\phi) = - \int_0^1 \int_M \ddt{\phi_t} (R_t-\underline{R})
\frac{\chi_{\phi_t}^n}{n!} dt,$$
where $\{ \phi_t \}_{0 \le t \le 1}$ is a path in $\mathcal{H}$ between $0$ and
$\phi$,
$R_t$ is the scalar curvature of $\chi_{\phi_t}$ and $\underline{R}$ is the
average of the scalar curvature.  The critical points of this functional are 
metrics of constant scalar curvature.

If $c_1(M)<0$ then there exists a 
K\"ahler-Einstein metric in the class $-c_1(M)$ (\cite{Y1}, \cite{Au}). 
It follows easily that the Mabuchi energy is bounded below in this class. 
Also, the
result of Donaldson \cite{D3} implies that 
$M$ is asymptotically Chow stable  with respect to the canonical
bundle.
  More generally, for any class, it is expected that the lower boundedness of the
Mabuchi energy should be equivalent to some notion of semistability in the sense
of geometric invariant theory (see
\cite{Y2}, \cite{T2,T3}, \cite{PS} and \cite{D4}).  

Chen \cite{C1} shows that if $c_1(M)$ is negative with $-\omega \in c_1(M)$ and
if there is a solution of (\ref{maineqn}) in $[\chi_0]$, then the Mabuchi energy
is  bounded below {\em in the class $[\chi_0]$}.
This suggests that (\ref{maineqn}) is related to how a change of
polarization affects the condition of stability of a manifold.
We immediately have the following corollary of the main theorem, generalizing
the result of Chen for dimension $n=2$. 

\bigskip
\noindent
{\bf Corollary.}
{\em Let $M$ be a compact K\"ahler manifold with $c_1(M)<0$.  If the
K\"ahler class $[\chi_0]$ satisfies the
inequality
\begin{equation} \label{classinequality}
- n \frac{c_1(M) \cdot [\chi_0]^{n-1}}{[\chi_0]^n} [\chi_0] + (n-1) c_1(M) >0,
\end{equation}
then the Mabuchi energy is bounded below on  $[\chi_0]$.}

\bigskip

Note that if $[\chi_0]=[-c_1(M)]$ then the inequality is more
than adequately satisfied, and so the set of $[\chi_0]$ satisfying the condition
is a reasonably large open set containing the canonical class.

In section 2, we prove a second order estimate of $\phi$ in terms of
$\phi$ itself, and in section 3, we prove the zero order estimate and complete
the proof of the main theorem.  The techniques used are generalizations of those
given in
\cite{W1}, and we will refer the reader to that paper for some of the calculations
and arguments.

In the following, $C_0, C_1, C_2, \ldots$ will denote constants depending only on
the initial data.

\renewcommand{\thesection}{2}
\addtocounter{section}{1}
\setcounter{equation}{0}
\bigskip
\bigskip
\noindent
{\bf 2. Second order estimate}
\bigskip

We use the maximum principle to prove the following estimate on the second
derivatives on $\phi$.

\begin{theorem} \label{theoremC2}
Suppose that
\begin{equation*}
nc \chi_0 - (n-1) \omega >0.
\end{equation*} 
Let $\phi=\phi_t$ be a solution of the J-flow (\ref{eqnJflow}) on $[0,
\infty)$.  Then there exist constants $A>0$ and $C>0$ depending only on the
initial data such that for any time $t\ge 0$, $\chi = \chi_{\phi_t}$ satisfies
\begin{equation} \label{eqnC2}
\Lambda_{\omega} \chi \le C e^{A (\phi - \inf_{M \times [0,t]} \phi)}.
\end{equation}
\end{theorem}
\addtocounter{lemma}{1}

\begin{proof}
We will assume that $\omega$ has been scaled so that $c=1/n$.  
Choose $0 < \epsilon < 1/(n+1)$ to be sufficiently small so that
\begin{equation} \label{eqnepsilon}
\chi_0 \ge (n-1 + (n+1)\epsilon) \omega. 
\end{equation}
We will use the
same notation as in \cite{W1}.  In particular, the operator $\lt$ acts on
functions $f$ by
$$\lt f = \frac{1}{n} \h{k}{l} \dd{k} \dbr{l} f, \qquad \textrm{where } \qquad
\h{k}{l} =
\chii{k}{j} \chii{i}{l} \g{i}{j}.$$
We calculate the evolution of $(\log(\Lambda_{\omega} \chi) - A
\phi)$, where $A$ is a constant to be determined.  From \cite{W1}, we have
\begin{eqnarray*}
\lefteqn{(\lt - \ddt{}) (\log (\Lambda_{\omega} \chi) - A\phi)} \\
& \ge & \mbox{} \frac{1}{n} (-C_0
\h{k}{l} \g{k}{l} - \frac{1}{\Lambda_{\omega} \chi} \chii{k}{l} R_{k
\overline{l}} - 2A \chii{i}{j} \g{i}{j} + A\h{k}{l} \chz{k}{l} + A) \\
& = & \frac{1}{n}(-C_0
\h{k}{l} \g{k}{l} - \frac{1}{\Lambda_{\omega} \chi} \chii{k}{l} R_{k
\overline{l}} - 2A \chii{i}{j} \g{i}{j} + (1-\epsilon) A\h{k}{l} \chz{k}{l} \\
&& \mbox{} +
\epsilon A \h{k}{l} \chz{k}{l} + A),
\end{eqnarray*}
where $C_0$ is a lower bound for the bisectional curvature of $\omega$, and
$R_{k \overline{l}}$ is the Ricci curvature tensor of $\omega$.  Recall from (2.4)
in
\cite{W1}, that
$\chi$ is uniformly bounded away from zero.  Hence we can choose $A$ to be large
enough so that
$$ \epsilon A \h{k}{l} \chz{k}{l} \ge C_0 \h{k}{l} \g{k}{l} +
\frac{1}{\Lambda_{\omega} \chi} \chii{k}{l} R_{k
\overline{l}}.$$
Now fix
a time $t>0$.  There is a point $(x_0, t_0)$ in $M \times [0,t]$ at which the
maximum of $(\log(\Lambda_{\omega} \chi) - A
\phi)$ is achieved.  We may assume that $t_0>0$.   Then at this point $(x_0,
t_0)$, we have
$$  1 + (1-\epsilon) \h{k}{l} \chz{k}{l} - 2 \chii{i}{j} \g{i}{j}\le 0.$$
From (\ref{eqnepsilon}), we get
$$ 1 + (n-1 + \epsilon) \h{k}{l} \g{k}{l} - 2\chii{i}{j} \g{i}{j} \le 0.$$
We will compute in normal coordinates at $x_0$ for $\omega$ in which $\chi$ is
diagonal and has eigenvalues $\lambda_1 \le \ldots \le \lambda_n$.  The above
inequality becomes
\begin{equation} \label{maininequality}
1 + (n-1 + \epsilon) \sum_{i=1}^n
\frac{1}{\lambda_i^2} - 2 \sum_{i=1}^{n} \frac{1}{\lambda_i} \le 0.
\end{equation}
We claim that (\ref{maininequality}) gives an upper bound for the
$\lambda_i$. 
To see this, complete the square to obtain
$$\sum_{i=1}^n \left(
\frac{1}{\sqrt{n-1+\epsilon}}  - \frac{\sqrt{n-1+\epsilon}}{\lambda_i}
\right)^2 \le \frac{n}{n-1+\epsilon} -1.$$ Hence, for $i=1, \ldots, n$,
$$  \frac{1}{\sqrt{n-1+\epsilon}} -\frac{\sqrt{n-1+\epsilon}}{\lambda_i}\le 
\frac{\sqrt{1-\epsilon}}{\sqrt{n-1+\epsilon}},$$
from which we obtain the upper bound
$$\lambda_i \le \frac{n-1+\epsilon}{1-\sqrt{1-\epsilon}}.$$
Hence at the point $(x_0, t_0)$, we have a constant $C$ depending only on
the initial data such that
$$\Lambda_{\omega} \chi \le C.$$
Then, on $M \times [0,t]$,
$$\log (\Lambda_{\omega} \chi) - A\phi \le \log C - A\inf_{M \times [0,t]}
\phi.$$
Exponentiating gives
$$\Lambda_{\omega} \chi \le C e^{A(\phi - \inf_{M \times [0,t]} \phi)},$$
completing the proof of the theorem.
\end{proof}

\renewcommand{\thesection}{3}
\addtocounter{section}{1}
\setcounter{theorem}{0}
\setcounter{lemma}{0}
\setcounter{equation}{0}
\bigskip
\bigskip
\noindent
{\bf 3. Proof of the Main Theorem}
\bigskip

We know from \cite{C2} that the flow exists for all time.  To prove the main
theorem we need uniform estimates on $\phi_t$ and all of its derivatives.
Given such estimates, the argument of section 5 of \cite{W1},
which is valid for any dimension, shows that $\phi_t$ converges in
$C^{\infty}$ to a smooth critical metric.

From Theorem \ref{theoremC2}, and standard parabolic methods, it suffices
to have a uniform $C^0$ estimate on $\phi$.  We prove this below, generalizing
the method of
\cite{W1}, using the precise form of the estimate (\ref{eqnC2}) and a Moser
iteration type argument.

\begin{theorem} \label{theoremC0}
\addtocounter{lemma}{1}
Suppose that
\begin{equation*} 
nc \chi_0 - (n-1) \omega >0.
\end{equation*} 
Let $\phi=\phi_t$ be a solution of the J-flow (\ref{eqnJflow}) on $[0,
\infty)$.  Then there exists a constant $\tilde{C}>0$ depending only on
the initial data such that 
$$\| \phi_t \|_{C^0} \le \tilde{C}.$$
\end{theorem}
\begin{proof}
We begin with a lemma.

\begin{lemma}  \label{lemmasup} $0 \le \sup_M \phi_t \le - C_1 \inf_M \phi_t +
C_2.$
\end{lemma}
\begin{proof} We will use the functional $I_{\chi_0}$ defined on
$\mathcal{H}$ by
\begin{equation} \label{eqnI}
I_{\chi_0}(\phi) = \int_0^1 \int_M \ddt{\phi_t} \frac{\chi_{\phi_t}^n}{n!}dt,
\end{equation}
for $\{\phi_t\}$ a path between $0$ and $\phi$ (this is a well-known functional,
see  \cite{D2} for example).   Taking the path
$\phi_t = t\phi$, we obtain the formula:
\begin{eqnarray} 
\nonumber
I_{\chi_0}(\phi) & = & \frac{1}{n!} \int_0^1 \int_M \phi \, \chi_{t \phi}^n dt \\
\nonumber
& = & \frac{1}{n!} \int_0^1 \int_M \phi \, (t \chi_{\phi} + (1-t) \chi_0)^n dt \\
\label{eqnI2}
& = & \frac{1}{n!} \sum_{k=0}^n \left[ \left(\genfrac{}{}{0pt}{}{n}{k}  \right)
\int_0^1 t^k(1-t)^{n-k}dt \right] \int_M \phi \, \chi_{\phi}^k \wedge
\chi_0^{n-k}.
\end{eqnarray}
From (\ref{eqnI}), we see that $I(\phi_t)=0$ along the flow.  The first
inequality then follows immediately, since the expression in
the square brackets in (\ref{eqnI2}) is a positive function of $n$ and $k$.  The
second inequality follows from (\ref{eqnI2}), the fact that $\triangle_{\omega}
\phi_t > - \Lambda_{\omega} \chi_0$, and properties of the Green's
function of $\omega$.
\end{proof}

From this lemma, it is sufficient to prove a lower bound for $\inf_M \phi_t$.  If
such a lower bound does not exist, then we can choose a sequence of times $t_i
\rightarrow \infty$ such that
\begin{enumerate}
\item[(i)] $\inf_M \phi_{t_i} = \inf_{t \in [0,t_i]} \inf_M \phi_t$
\item[(ii)] $\inf_M \phi_{t_i} \rightarrow -\infty.$
\end{enumerate}
We will find a contradiction.  Set $B = A/(1-\delta)$ where $A$ is the constant
from (\ref{eqnC2}), and let $\delta$ be a small positive constant to be determined
later.  Let 
$$\psi_{i} = \phi_{t_i} - \sup_M \phi_{t_i},$$
and let $u= e^{-B\psi_{i}}$.  We will show that $u$ is uniformly bounded
from above, which will give the contradiction.  First, we have the following
lemma.

\begin{lemma} For any $p \ge 1$,
\begin{equation} \label{estimateu}
\int_M | \nabla u^{p/2} |^2 \frac{\omega^n}{n!} \le C_3 \, p \, \|u
\|_{C^0}^{1-\delta}
\int_M u^{p-(1-\delta)} \frac{\omega^n}{n!}.
\end{equation}
\end{lemma}

\begin{proof} The proof is given for $n=2$ in \cite{W1}, and since the same argument
works for any dimension, we will not reproduce it here.  Crucially, the proof uses the estimate
(\ref{eqnC2}).
\end{proof}

We will use the notation
$$\| f \|_{c} = \left( \int_M |f|^c \frac{\omega^n}{n!} \right)^{1/c},$$
for $c>0$.  It is not a norm for $0 < c <1$ but this fact is not important.  The
following lemma allows us to estimate the
$C^0$ norm of
$u$ using a Moser iteration type method (compare to \cite{Y1}).

\begin{lemma}
If $u \ge 0$ satisfies the estimate (\ref{estimateu}) for all $p \ge 1$, then for
some constant $C'$ independent of $u$,
$$ \| u \|_{C^0} \le C' \| u \|_{\delta}.$$
\end{lemma}
\begin{proof}
For $\beta = n/(n-1)$, the Sobolev inequality for functions $f$ on $(M, \omega)$
is
$$\| f \|_{2 \beta}^2 \le C_4 (\| \nabla f \|_2^2 + \|f \|_2^2).$$
Applying this to $u^{p/2}$ and making use of (\ref{estimateu}) gives
$$ \|u \|_{p\beta} \le C_5^{1/p} p^{1/p} \|u \|_{C^0}^{\gamma/p} \| u
\|_{p-\gamma}^{(p - \gamma)/p},$$
for $\gamma = 1 - \delta$.  By replacing $p$ with $p \beta + \gamma$ we obtain
inductively
$$ \|u \|_{p_k \beta} \le C(k) \| u \|_{C^0}^{1-a(k)} \| u
\|_{p-\gamma}^{a(k)},$$
where
\begin{eqnarray*}
p_k & = & p \beta^k + \gamma (1 + \beta + \beta^2 + \cdots + \beta^{k-1}) \\
C(k) & = & C_5^{(1+ \beta + \cdots + \beta^k)/p_k} p_0^{\beta^k/p_k}
p_1^{\beta^{k-1}/p_k} \cdots p_k^{1/p_k} \\
a(k) & = & \frac{(p-\gamma)\beta^k}{p_k}.
\end{eqnarray*}
Set $p=1$.  Note that for some fixed $l$, $\beta^k \le p_k \le
\beta^{k+l}$.  It is easy to check that $C(k)
\le C_6$ for some constant $C_6$.  As $k$ tends to infinity, $p_k\rightarrow
\infty$, $a(k) \rightarrow a \in (0,1)$, and the required estimate follows
immediately. 
\end{proof}

We can now finish the proof of Theorem \ref{theoremC0}.
Since $u=e^{-B\psi_i}$ and $\psi_i$ satisfies $\sup_M \psi_i =0$ and
$$\chz{k}{l} + \dd{k} \dbr{l}\psi_i \ge 0,$$
we can apply Proposition 2.1 of \cite{T1} to get a bound on $\| u \|_{\delta}$
for $\delta$ small enough.  This completes the proof of
 Theorem \ref{theoremC0}.
\end{proof}

\bigskip
\noindent
{\bf Acknowledgements.} The author would like to thank his thesis advisor, D.H.
Phong for his constant support, advice and encouragement.   The author also
thanks Jacob Sturm, Lijing Wang, Mao-Pei Tsui and Lei Ni for some helpful
discussions.  In addition, he is very grateful to S.-T. Yau for his
support and advice.   The results of this paper are contained in the
author's PhD thesis at Columbia University
\cite{W2}.

\end{document}